\documentclass[twoside]{article}
\title{Weak Gorenstein global
dimension}
\date{}
\usepackage[all]{xy}
\usepackage{amsfonts}
\usepackage{amsmath}
\usepackage{amssymb}
\usepackage{latexsym}
\usepackage{mathrsfs}
\usepackage{eufrak}

\newtheorem{thm}{\bf Theorem}[section]
\newtheorem{cor}[thm]{\bf Corollary}
\newtheorem{lem}[thm]{\bf Lemma}
\newtheorem{prop}[thm]{\bf Proposition}
\newtheorem{defn}[thm]{\bf Definition}
\newtheorem{rem}[thm]{\bf Remark}
\newtheorem{exmp}[thm]{\bf Example}
\newtheorem{conj}[thm]{\bf Conjecture}
\catcode`\ç=13
\defç{\c{c}}
\catcode`\é=13
\defé{\'e}
\catcode`\à=13
\defà{\`a}
\catcode`\è=13
\defè{\`e}
\catcode`\â=13
\defâ{\^a}
\catcode`\ù=13
\defù{\`u}
\catcode`\ê=13
\defê{\^e}
\catcode`\î=13
\defî{\^\i}
\catcode`\ô=13
\defô{\^o}

\def\proof{{\parindent0pt {\bf Proof.\ }}}

\def\fd{{\rm fd}}
\def\id{{\rm id}}

\def\Gfd{{\rm Gfd}}

\def\Im{{\rm Im}}
\def\Coker{{\rm Coker}}
\def\Ker{{\rm Ker}}

\def\Ext{{\rm Ext}}
\def\Tor{{\rm Tor}}
\def\Hom{{\rm Hom}}

\def\sup{{\rm sup}}

\def\wdim{{\rm wdim}}

\def\lGwdim{{l\rm.Gwdim}}
\def\rGwdim{{r\rm.Gwdim}}

\def\lFFD{{l\rm.FFD}}
\def\lFGFD{{l\rm.FGFD}}
\def\rFFD{{r\rm.FFD}}
\def\rFGFD{{r\rm.FGFD}}

\newcommand{\cqfd}
{\hspace{1cm}
\rule{2mm}{2mm}%
\medbreak%
\par%
}

\begin{document}
\thispagestyle{empty}

\maketitle \vspace*{-1.5cm}

\begin{center}
{\large\bf Driss Bennis}

\bigskip

Department of Mathematics, Faculty of Science and Technology of
Fez,\\ Box 2202, University S. M. Ben Abdellah Fez,
Morocco, \\[0.2cm]
 driss\_bennis@hotmail.com
\end{center}

\bigskip\bigskip
\noindent{\large\bf Abstract.} In this paper, we investigate the
weak Gorenstein global dimensions.  We are mainly interested in
studying the problem when the left and right weak Gorenstein
global dimensions coincide. We first show, for GF-closed rings,
that the left and right weak Gorenstein global dimensions are
equal when they are finite. Then, we prove that the same equality
holds for any two-sided coherent ring.  We conclude with some
examples and a brief discussion of the scope and limits of our
results\bigskip

\begin{section}{Introduction} Throughout the paper, all rings
are associative with identity, and all modules are unitary.\\
Let $R$ be a ring. The injective (resp., flat) dimension of an
$R$-module $M$ is denoted by $\id_R(M)$ (resp.,
$\fd_R(M)$).\\
\indent  A left (resp., right) $R$-module  $M$ is called
\textit{Gorenstein flat}, if there exists an exact sequence of
flat left (resp., right) $R$-modules  $$\mathrm{F}=\
\cdots\rightarrow F_1\rightarrow F_0 \rightarrow F^0 \rightarrow
F^1 \rightarrow\cdots ,$$ such that $M \cong \Im(F_0 \rightarrow
F^0)$ and such that the sequence $I\otimes_R \mathrm{F} $
 (resp.,  $\mathrm{F}\otimes_R I $) remains exact whenever
$I$ is an injective  right (resp., left) $R$-module. The sequence $\mathrm{F}$ is called a complete flat resolution.\\
\indent For a positive integer $n$, we say that $M$ has
\textit{Gorenstein flat dimension} at most $n$, and we write
$\Gfd_R(M)\leq n$, if there is an exact sequence of $R$-modules $$
0 \rightarrow G_n\rightarrow \cdots \rightarrow G_0\rightarrow M
\rightarrow 0,$$ where each $G_i$ is Gorenstein flat (please see
\cite{LW,Rel-hom,HH}).\bigskip

The notion of Gorenstein flat modules was introduced and studied
over Gorenstein rings, by  Enochs, Jenda, and Torrecillas
\cite{GoPlat}, as a generalization of the notion of flat modules
in the sense that an $R$-module is flat if  and only if  it is
Gorenstein flat with finite flat dimension. In \cite{FCring}, Chen
and Ding generalized known characterizations of Gorenstein flat
modules (then of the Gorenstein flat dimension) over Gorenstein
rings to $n$-FC rings (coherent with finite self-FP-injective
dimension). Then, in \cite{HH}, Holm generalized the study of the
Gorenstein flat dimension to coherent rings. In the same
direction,     the study of Gorenstein flat dimension is
generalized, in \cite{B}, to a larger class of rings called
GF-closed: a ring $R$ is called left (resp., right) GF-closed, if
the class of all Gorenstein flat left (resp., right)
 $R$-modules is closed under extensions; that is, for every short
exact sequence of left (resp., right) $R$-modules $0\rightarrow A
    \rightarrow B \rightarrow C \rightarrow 0$, the condition $A$
and $C$ are Gorenstein flat implies that $B$ is   Gorenstein flat.
A ring is called GF-closed, if it is both left and right
GF-closed. The class of GF-closed rings includes strictly the one
of coherent rings and also the one of rings of finite weak global
dimension \cite[Example 3.6]{B}.\bigskip

In this paper, we are concerned with the left and right weak
Gorenstein global dimensions of rings, $\lGwdim(R)$ and
$\rGwdim(R)$,  which are respectively defined as follows:\medskip

\noindent  \mbox{}\hspace{2.5cm}    $\lGwdim(R) =
\sup\{\Gfd_R(M)\,|\,M\ is\ a\ left\ R\!-\!module\}\ $  and\\
 \mbox{}\hspace{2.5cm}
 $\rGwdim(R) = \sup\{\Gfd_R(M)\,|\,M\ is\  a\ right \
R\!-\!module\}.$\bigskip

In the classical case we have, for any ring $R$, the following
well-known equality \cite[Theorem 9.15]{Rot}:
 $$\sup\{\fd_R(M)\,|\,M\ is\  a\ right \
R\!-\!module\}=\sup\{\fd_R(M)\,|\,M\ is\ a\ left\
R\!-\!module\}.$$ The common value of these equal terms is called
weak global dimension
of $R$ and denoted by $\wdim(R)$.\\
In \cite[Theorem 2.2]{B1}, it is proved, for any ring $R$, that
$\Gfd_R(M)= \fd_R(M)$ for every  (left or  right) $R$-module $M$
with finite flat dimension. Then, if $\wdim(R)<\infty$, we get:
$$\lGwdim(R) = \wdim(R) =\rGwdim(R).$$ On the
other hand, by \cite[Theorem 12.3.1 ($1 \Leftrightarrow
4$)]{Rel-hom}, we have for a left and right Noetherian ring $R$:
$$\lGwdim(R) = \rGwdim(R).$$
This naturally leads to the following conjecture:

\begin{conj}\label{conje} For any ring $R$, $\lGwdim(R) = \rGwdim(R)$.
\end{conj}

The main purpose of this paper is to prove that this conjecture
holds for a large class of rings. First, we prove that the
conjecture is true for GF-closed rings which have finite both left
and right weak Gorenstein global dimensions (Theorem \ref{thm-1}).
Then, we prove the conjecture is true for two-sided coherent rings
(Theorem \ref{thm-2}). In Proposition \ref{prop-IF}, we prove
that, for a ring $R$, $\lGwdim(R) =0$ if and only if
$\rGwdim(R)=0$; and, in this case, $R$ is an IF ring (i.e., $R$
satisfies:  every injective right (resp., left) $R$-module $I$ is
flat).  We conclude with some examples  and a brief discussion of
the scope and limits of our results (Remark \ref{rem}, Proposition
\ref{prop-product}, and Example \ref{exm}).

\end{section}


\begin{section}{Main results}

We begin with the first main result which says  that Conjecture
\ref{conje} is true for GF-closed rings with  finite both left and
right weak Gorenstein global dimensions.

\begin{thm}\label{thm-1}
If $R$ is a   GF-closed ring with    finite both left and right
weak Gorenstein global dimensions, then $\lGwdim(R)= \rGwdim(R).$
\end{thm}

To prove this theorem, we need the following results.\\
\indent The following lemma  generalizes \cite[Lemma 2.19]{CFH}.

\begin{lem}\label{lem-1} Assume that $R$ is a  left (resp., right) GF-closed  ring. If
$M$ is a left  (resp., right) $R$-module with $\Gfd_R(M)<\infty$,
then there exists a short exact sequence of left  (resp., right)
$R$-modules $0\rightarrow M\rightarrow M'\rightarrow G \rightarrow
0$, such that $\fd_R(M')=\Gfd_R(M)$ and $G$ is Gorenstein flat.
\end{lem}
\proof We only prove the case of left modules, and the  case of
right modules is proved similarly.\\
Let $\Gfd_R(M)=n$ for some positive integer $n$. We prove the result
by induction on $n$. The case $n=0$ holds by the definition of the
Gorenstein flat module. Then, suppose that $n>0$ and pick  a short
exact sequence of left $R$-modules:
 $0\rightarrow K\rightarrow G_0\rightarrow M \rightarrow 0,$ where
$G_0$ is Gorenstein flat and $\Gfd_R(K)=n-1$. By induction, there
exists a short exact sequence of $R$-modules: $0\rightarrow
K\rightarrow K'\rightarrow H \rightarrow 0,$  such that
$\fd_R(K')=\Gfd_R(K)=n-1$ and $H$ is Gorenstein flat. Consider the
pushout   diagram:
$$\xymatrix{
     &   0 \ar[d] & 0 \ar[d]  &   &  \\
0\ar[r]& K \ar[d] \ar[r] & G_{0} \ar@{-->}[d] \ar[r] & M \ar@{=}[d]
\ar[r] &
0\\
0\ar[r]& K' \ar[d]  \ar@{-->}[r] &D\ar[d] \ar[r] & M
\ar[r] & 0\\
 &  H \ar[d] \ar@{=}[r] &H \ar[d]  &   &  \\
     &   0 & 0   &   & }$$
By the middle vertical sequence and since $R$ is  left GF-closed,
$D$ is Gorenstein flat. Then, there exists  a short exact sequence
of left $R$-modules $0\rightarrow D\rightarrow F\rightarrow G
\rightarrow 0,$ where $F$ is flat and $G$ is Gorenstein flat.
Consider  the pushout diagram:
$$\xymatrix{
     &    &  0 \ar[d] & 0 \ar[d]  &  \\
0\ar[r]& K' \ar@{=}[d] \ar[r] &D\ar[d] \ar[r] & M \ar@{-->}[d]
\ar[r] &
0\\
0\ar[r]& K'   \ar[r] & F\ar[d] \ar@{-->}[r] & M' \ar[d] \ar[r] &
0\\
 &   & G\ar[d] \ar@{=}[r] & G\ar[d]   &
 \\
 &   & 0 &0  & }$$ By the middle  horizontal sequence $\fd_R(M')=n$.
Then, the right vertical sequence is the desired
sequence.\cqfd\bigskip

Compare the following result to \cite[Theorem 2.6 (ii)]{HH1}.
\begin{cor}\label{cor-1} Assume that $R$ is a  left (resp., right) GF-closed  ring. If
$M$ is an  injective left  (resp., right) $R$-module, then $\fd_R(M
)=\Gfd_R(M)$.
\end{cor}
\proof It is known that  $\Gfd_R(M )\leq \fd_R(M)$ for every (left
or right) $R$-module and over any associative ring $R$. Conversely,
assume that $\Gfd_R(M)$ is finite, then by Lemma \ref{lem-1} there
exists a short exact sequence of  $R$-modules:
 $0\rightarrow M\rightarrow M'\rightarrow G \rightarrow 0$ such
that $\fd_R(M')=\Gfd_R(M)$. Since $M$ is injective, this sequence
splits and therefore $\fd_R(M)\leq \fd_R(M')=\Gfd_R(M)$.\cqfd

\begin{lem}\label{lem-2}
If $R$ is a left GF-closed ring with $\lGwdim(R)<\infty$, then, for
a positive integer $n$, the following are equivalent:
\begin{enumerate}
    \item $\lGwdim(R)\leq n$;
    \item $\Gfd_{R}(M)\leq n$ for every finitely presented left $R$-module $M$;
     \item $\Gfd_{R}(R/I)\leq n$ for every   finitely generated left ideal $I$ of $R$;
    \item   $ \fd_R(E)\leq n $ for every injective right $R$-module $E$;
   \item   $ \fd_R(E')\leq n $ for every right $R$-module $E'$ with finite injective dimension;
\end{enumerate}
Consequently, the left weak Gorenstein  global dimension of $R$ is
also determined by the formulas:
\begin{eqnarray}
 \nonumber  \lGwdim(R)&=& \sup\{\Gfd_{R}(R/I)\,|\,I\ is \ a\  finitely\ generated\ left \ ideal\ of\ R\} \\
 \nonumber          &=& \sup\{\Gfd_{R}(M)\,|\,M\ is \ a\ finitely\ presented\ left \ R\!-\!module \} \\
\nonumber         &=& \sup\{\fd_{R}(E)\,|\,E\ is \ an\ injective\ right\ R\!-\!module \} \\
 \nonumber        &=& \sup\{\fd_{R}(E')\,|\,E'\ is \ a\ right\ R\!-\!module \
 with\
 \id_R(E)<\infty\}.
\end{eqnarray}
\end{lem}
\proof The implications $1\Rightarrow 2\Rightarrow 3$ are trivial.
The implication $3\Rightarrow 4$ follows from \cite[Theorem 2.8
($1\Rightarrow 2$)]{B}. The implication $4\Rightarrow 5$ is proved
by induction on $ \id_R(E')$ using the flat counterpart of
\cite[Corollary 2, p. 135]{Bou}. Finally, the implication
$5\Rightarrow 1$ is a simple consequence of \cite[Theorem 2.8
($3\Rightarrow 1$)]{B}.\cqfd\bigskip

Similarly we obtain the right version of Lemma \ref{lem-2}.

\begin{lem}\label{lem-3}
If $R$ is a right GF-closed ring with $\rGwdim(R)<\infty$, then, for
a positive integer $n$, the following are equivalent:
\begin{enumerate}
    \item $\rGwdim(R)\leq n$;
    \item $\Gfd_{R}(M)\leq n$ for every finitely presented right $R$-module $M$;
     \item $\Gfd_{R}(R/I)\leq n$ for every   finitely generated right  ideal $I$ of $R$;
    \item   $ \fd_R(E)\leq n $ for every injective left $R$-module $E$;
   \item   $ \fd_R(E')\leq n $ for every left $R$-module $E'$ with finite injective dimension;
\end{enumerate}
Consequently, the right weak Gorenstein global dimension of $R$ is
also determined by the formulas:
\begin{eqnarray}
 \nonumber  \rGwdim(R)&=& \sup\{\Gfd_{R}(R/I)\,|\,I\ is \ a\ finitely\ generated\ right \ ideal\ of\ R\} \\
 \nonumber          &=& \sup\{\Gfd_{R}(M)\,|\,M\ is \ a\ finitely\ presented\ right \ R\!-\!module \} \\
\nonumber         &=& \sup\{\fd_{R}(E)\,|\,E\ is \ an\ injective\ left\ R\!-\!module \} \\
 \nonumber        &=& \sup\{\fd_{R}(E')\,|\,E'\ is \ a\ left\ R\!-\!module\ with\ \id_R(E)<\infty\}.\\
 \nonumber
\end{eqnarray}
\end{lem}

\noindent\textbf{Proof of Theorem \ref{thm-1}.} Assume that
$\rGwdim(R)=n$ is finite. From Corollary \ref{cor-1}, $\fd_R(E
)=\Gfd_R(E)\leq n$ for every injective right $R$-module  $E$.
Then, from
Lemma \ref{lem-2}, $\lGwdim(R)\leq n= \rGwdim(R).$\\
The converse inequality is proved similarly.\cqfd\bigskip

Under the condition of Theorem \ref{thm-1}, the classical left and
right finitistic flat dimension are equal and they are also equal
to the left and right weak Gorenstein global dimensions. Recall
that the left finitistic flat dimension, $\lFFD(R)$, of a ring $R$
is defined as follows: $$\lFFD(R)= \{ \fd_R(M)\, | \, M\ is \ a\
left\ R\!-\!module\ with \ \fd_R(M)<\infty\}.$$ The right
finitistic dimension  $\rFFD(R)$ is defined similarly.

\begin{prop}\label{prop-Gwgldim-FFD}
If $R$ is a  GF-closed ring with      finite both left and right
weak Gorenstein global dimensions, then $\lFFD(R)=\lGwdim(R)=
\rGwdim(R)=\rFFD(R).$
\end{prop}
\proof This follows from Theorem \ref{thm-1} and the following
result.\cqfd\bigskip

Recall that the left finitistic Gorenstein flat dimension,
$\lFGFD(R)$, of   a ring $R$ is defined as follows: $$\lFGFD(R)=
\{ \Gfd_R(M)\, | \, M\ is \ a\ left\ R\!-\!module\ with \
\Gfd_R(M)<\infty   \}.$$ The right  finitistic dimension  $\rFGFD(R)$ is defined similarly.\\
The following result  is a generalization of \cite[Theorem
3.24]{HH}.
\begin{prop}\label{prop-FGFD} For any ring $R$, we have  $\lFFD(R)\leq\lFGFD(R)$ and $\rFFD(R)\leq\rFGFD(R)$.\\
Furthermore, if $R$ is  left (resp., right) GF-closed, then
$\lFFD(R)=\lFGFD(R)\ (resp., \ \rFFD(R)=\rFGFD(R)).$
\end{prop}
\proof The inequalities follow immediately by the fact that
$\Gfd_R(M)= \fd_R(M)$ for every
$R$-module $M$ with finite flat dimension \cite[Theorem 2.2]{B1}.\\
Now, assume that $R$ is  left GF-closed (the right version is proved
similarly). It remains to prove the converse inequality
$\lFGFD(R)\leq\lFFD(R)$. For that, we can assume that $\lFFD(R)=n$
is finite. Let $M$ be a left $R$-module with finite Gorenstein flat
dimension. By Lemma \ref{lem-1}, there exists a short exact sequence
of left $R$-modules $0\rightarrow M\rightarrow M'\rightarrow G
\rightarrow 0$ such that $\Gfd_R(M)=\fd_R(M')\leq n$. This implies
the desired inequality.\cqfd\bigskip

Now we give the second main result which says  that Conjecture
\ref{conje} is true for two-sided coherent rings. For that we use the following notions:\\
From \cite{Colby}, a ring $R$ is called right (resp., left) IF, if
every injective right (resp., left) $R$-module $I$ is flat. A ring
$R$ is called IF, if it is both left and right IF. Then, let us
call a ring $R$ is right (resp., left) $n$-IF, for $n\geq 0$, if $
\fd_R(E)\leq n $ for every injective right (resp., left)
$R$-module $E$. And $R$ is called $n$-IF, if it is both left and
right
$n$-IF.\\
Obviously, $0$-IF rings are just the IF rings. And, from
\cite[Theorem 9.1.11]{Rel-hom}, the $n$-IF Noetherian rings are
the
same as the well-known $n$-Gorenstein rings.\\
The following result is also a generalization of \cite[Theorem
12.3.1 ($1 \Leftrightarrow 4$)]{Rel-hom}.

\begin{thm}\label{thm-2}
If $R$ is a right and left coherent ring, then, for a positive
integer $n$, the following are equivalent:
\begin{enumerate}
    \item $\lGwdim(R)\leq n$;
        \item   $R$ is $n$-IF;
    \item $\rGwdim(R)\leq n$.
\end{enumerate}
Consequently, for any two-sided coherent ring $R$, $\lGwdim(R)=
\rGwdim(R)$.
\end{thm}

The proof of this theorem uses the notion of a flat preenvelope of
modules which is defined as follows:

\begin{defn}[\cite{Xu}]\label{def-prenvel}
\textnormal{ Let $R$ be a ring and let $F$ be a flat $R$-module. For
an $R$-module $M$, an homomorphism  (or $F$) $\varphi :\,
M\rightarrow F$ is called  a flat preenvelope, if for any
homomorphism $\varphi' :\, M\rightarrow F'$ with $F'$ is a flat
module, there is an homomorphism $f :\, F\rightarrow F'$ such that
$\varphi' =f\varphi$.}
\end{defn}

Note that if $M$ embeds in a flat module, then its flat preenvelope
(if it exists) is injective.\\
\indent The coherent rings is also characterized by the notion of
a flat preenvelope of modules as follows:

\begin{lem}[\cite{Xu}, Theorem  2.5.1]\label{car-coherent-prevelo}
A ring $R$ is coherent if and only if every $R$-module has a flat
preenvelope.
\end{lem}

Also we use the notion of a flat cover of modules which is defined
as follows:

\begin{defn}[\cite{Xu}]\label{def-cov}
\textnormal{Let $R$ be a ring and let $F$ be a flat $R$-module. For
an $R$-module $M$, an homomorphism (or $F$) $\varphi :\,
F\rightarrow M$ is called  a flat precover, if for any homomorphism
$\varphi' :\, M\rightarrow F'$ with $F'$ is a flat module, there is
an homomorphism $f :\, F'\rightarrow F $ such that $\varphi'
=\varphi
f$.\\
A flat precover $\varphi :\, F\rightarrow M$ of $M$ is called flat
 cover, if   every endomorphism $f$  of $F$  with  $\varphi  =\varphi
 f$ must be an automorphism.}
\end{defn}

Recall that an $R$-module $M$ is called cotorsion, if
$\Ext^{1}_{R}(F, M)=0$ for every  flat $R$-module $F$.

\begin{lem}[\cite{BBE} and \cite{Rel-hom}, Lemma 5.3.25]\label{lem-fCovers}
For any ring $R$, every $R$-module $M$ has a flat cover $\varphi :\,
F\rightarrow M$ such that $\Ker(\varphi)$ is cotorsion.
\end{lem}

Note that every flat cover is surjective.\\
\indent From its proof, \cite[Proposition 3.22]{HH} is stated, as
we need here, as follows:

\begin{lem}\label{lem-G-cot-fd}
Let $R$ be  a right (resp., left) coherent ring. If $T$ is a left
(resp., right) $R$-module such that $\Tor_1^R(I,T)=0$ (resp.
$\Tor_1^R(T,I)=0$) for every injective right (resp., left)
$R$-module, then $\Ext_R^1(T,K)=0$  for every cotorsion left
(resp., right) $R$-module $K$ with finite flat dimension.
\end{lem}

\noindent\textbf{Proof of Theorem \ref{thm-2}.} Since every coherent
ring is GF-losed, the implications $1\Rightarrow 2$ and
$3\Rightarrow 2$ follow from Theorem
\ref{thm-1} and Lemmas \ref{lem-2} and \ref{lem-3}.\\
We only prove the implication $2\Rightarrow 1$. The    implication
$3\Rightarrow 1$ has a similar proof.\\
Let $M$ be a left $R$-module, and consider an exact sequence of left
$R$-modules:
$$(*)\quad 0\rightarrow G \rightarrow P_{n-1}\rightarrow\cdots
\rightarrow P_0\rightarrow M \rightarrow 0,$$ where each $P_i$ is
projective. We have to prove that $G$ is Gorenstein flat.\\
First note that, using the above sequence $(*)$, we have:
$$\Tor^R_{k}(E,G)\cong \Tor^R_{n+k}(E,M
)\qquad  \mathrm{for}\ \mathrm{every} \ k\geq 1\ \mathrm{and}\
\mathrm{every}\   \mathrm{right}\ \mathrm{R \frac{\;}{}module}\
E.$$ If $E$ is an injective right $R$-module, then  $ \fd_R(E)\leq
n $ (since $R$ is $n$-IF), and so by the above isomorphism we get:
$$(**)\quad\Tor^R_{k}(E,G) =0\qquad  \mathrm{for}\ \mathrm{every} \
k\geq 1\ \mathrm{and}\ \mathrm{every}\ \mathrm{injective}\
\mathrm{right}\ \mathrm{R\frac{\;}{}module}\ E.$$ Then, by
\cite[Theorem 3.6 ($i \Leftrightarrow iii$)]{HH}, it remains to
construct a right flat resolution of $G$:
$$\mathbf{F}= 0\rightarrow G \rightarrow F^0\rightarrow
F^1\rightarrow\cdots,
$$ such that the sequence $\Hom_R( \mathbf{F}, F )$ is exact
whenever $F$ is a flat left $R$-module. Equivalently, for every
positive integer $i $, $G^i \rightarrow F^{i}$ is a flat
preenvelope of $G^i $, where $G^0=G$
and $G^i=\Ker(F^i \rightarrow F^{i+1})$ for $i\geq 1$.\\
Consider a short exact sequence of left $R$-modules $0 \rightarrow G
\rightarrow I \rightarrow L\rightarrow 0$, where $I$ is injective.
From Lemma \ref{lem-fCovers}, there exists a short exact sequence of
left $R$-modules  $0 \rightarrow K\rightarrow F \rightarrow
I\rightarrow 0$, where $F$ is flat and $K$ is cotorsion. Since $R$
is $n$-IF, $ \fd_R(I)\leq n $ and so $\fd_R(K)\leq n-1 $. Consider
the pullback diagram
$$\xymatrix{
     &  0 \ar[d]  & 0 \ar[d]  &   &  \\
 &  K\ar[d] \ar@{=}[r]&  K\ar[d]  &   &  \\
 0\ar[r]& D\ar@{-->}[d] \ar@{-->}[r] & F\ar[d] \ar[r] & L \ar@{=}[d]  \ar[r] & 0\\
0\ar[r]&G \ar[r] \ar[d]& I\ar[d] \ar[r]& L\ar[r] & 0\\
 & 0 &0  &   & }$$
From Lemma \ref{lem-G-cot-fd} and $(**)$, we have $\Ext_R^1(G,K)=0$.
Then, the left vertical sequence  splits, and so $G$ embeds in the
flat module $F$. Thus, $G$ admits an injective flat preenvelope
$G\rightarrow F^0$, which gives the desired first  flat
preenvelope.\\
Now, for $G^1=\Coker(G\rightarrow F^0)$ we prove that
$\Ext_R^1(G^1,K)=0$   for every cotorsion left $R$-module $K$ with
finite flat dimension. This gives, using the same argument above,
the desired second flat preenvelope $G^1\rightarrow F^1$, and
recursively we obtain the remains flat
preenvelopes.\\
Let $K$  be a cotorsion left $R$-module with finite flat dimension.
By Lemma \ref{lem-fCovers}, there exists a flat cover  $F
\rightarrow K$ of $K$ such that we obtain a short exact sequence  of
left $R$-modules   $$0 \rightarrow K'\rightarrow F \rightarrow
K\rightarrow 0,$$ where $K'$ is cotorsion with finite flat dimension
(since $\fd_R(K)<\infty$). By \cite[Proposition 3.1.2]{Xu}, $F$ is
cotorsion. Then, we get the following commutative diagram

$$\xymatrix{
    & 0 \ar[d]  &  0 \ar[d]  & 0 \ar[d]  &   \\
 0\ar[r]&  \Hom(G^1,K')\ar[d]\ar[r]  & \Hom(F^0,K')\ar[d]\ar[r] &  \Hom(G,K')\ar[d]  &   \\
 0\ar[r]& \Hom(G^1,F)\ar[d]\ar[r]  & \Hom(F^0,F)\ar[d]\ar[r]  & \Hom(G,F)\ar[d]\ar[r]  & 0\\
 0\ar[r]& \Hom(G^1,K ) \ar[r]  & \Hom(F^0,K )\ar[d]\ar[r]  & \Hom(G,K )\ar[d]   &  \\
 &  &0  &  0 & }$$
with exact rows and columns. Indeed, the middle vertical sequence
is exact since $F \rightarrow K$ is a flat cover of $K$; the right
vertical sequence is exact since $\Ext_R^1(G,K')=0$  since $K'$ is
cotorsion  with finite flat dimension; and the middle horizontal
sequence is exact since $G \rightarrow F^0$ is a flat preenvelope.
Then,   the   sequence $ 0 \rightarrow \Hom(G^1,K
)\rightarrow\Hom(F^0,K )\rightarrow\Hom(G,K )\rightarrow0$ is
exact. This implies, using $\Ext_R^1(F^0,K )=0$, that
$\Ext_R^1(G^1,K)=0$, and this completes the proof.\cqfd\bigskip

For the case where $\lGwdim(R) =0$ or $\rGwdim(R)=0$ we have the
following generalization of \cite[Theorem 6 ($1 \Leftrightarrow
2$)]{FCring}.

\begin{prop}\label{prop-IF}
For a ring  $R$, the following are equivalent:
\begin{enumerate}
    \item $\lGwdim(R) =0$;
        \item   $R$ is IF;
    \item $\rGwdim(R)=0$.
\end{enumerate}
\end{prop}
\proof the implications $1\Rightarrow 2$ and $3\Rightarrow 2$ follow
from Theorem \ref{thm-2}.\\
We prove the implication $2\Rightarrow 1$. The implication
$3\Rightarrow 1$ has a similar proof.\\
We prove that every left $R$-module $M$ is Gorenstein flat. For
that, we have to construct a complete flat resolution $\mathbf{F}$
such that $M\cong \Im(F_0 \rightarrow F^0)$. Since $R$ is IF, we can
consider any flat resolution of $M$ as the ``left half" of
$\mathbf{F}$. For that ``right half" of $\mathbf{F}$, consider an
injective resolution $\mathbf{I}$ of $M$. Since $R$ is IF, the
sequence  $\mathbf{I}$ is a right flat resolution of $M$ such that
the sequence $I\otimes_R \mathbf{I}$ is exact whenever $I$ is an
injective (then flat) right $R$-module, as desired.\cqfd

\begin{rem}\label{rem}\textnormal{
In \cite[Example 2]{Colby}, Colby gave an example of a left and
right coherent ring $R$ which is right IF but not left IF. Then,
$\lGwdim(R) =\rGwdim(R) = \infty $. Indeed, if
$\lGwdim(R)<\infty$, then, by Lemma \ref{lem-2}, $\lGwdim(R)=0$,
and from Theorem \ref{thm-2}, $\lGwdim(R) =\rGwdim(R)=0$. But,
this contradicts the fact that $R$ is not  left IF (Proposition
\ref{prop-IF}). Consequently:
\begin{enumerate}
    \item To have the implication $4 (\textnormal{or}\ 5) \Rightarrow 1$ in  Lemma \ref{lem-2}, the condition
    ``$\lGwdim(R)<\infty$" can not be dropped.
    \item In \cite[Theorem 3.14]{HH}, the condition
    ``$\Gfd(M)<\infty$"
can not be dropped to get the implication $iii (\textnormal{or}\ ii)
\Rightarrow i$. Indeed, since $\lGwdim(R)=\infty$, there exists,
using \cite[Proposition 3.13]{HH}, a left $R$-module with
$\Gfd_R(M)=\infty$. However, $\Tor^R_{i}(I,M)=0$ for every $i>0$ and
every injective (then flat) right $R$-module $I$.
\end{enumerate}}
\end{rem}

Finally, to construct examples of GF-closed rings of finite weak
Gorenstein global dimension which are neither coherent nor of
finite weak global dimension, we   need the following result:

\begin{prop}\label{prop-product}
For any family of  rings $\{R_{i}\}_{i=1,...,m}$, we have: \\
 \mbox{}\hspace{2.5cm}
 $\lGwdim(\displaystyle\prod_{i=1}^m R_i)=\sup\{\lGwdim(R_i), 1\leq i
\leq m\} \  $ and\\
 \mbox{}\hspace{2.5cm}    $\rGwdim(\displaystyle\prod_{i=1}^m
R_i)=\sup\{\rGwdim(R_i), 1\leq i \leq m\}.$
\end{prop}
\proof The result is a   consequence of \cite[Theorem
3.4]{B}.\cqfd

\begin{exmp}\label{exm}
\textnormal{Consider an IF ring $R_1$ with infinite weak global
dimension, and consider a non-coherent ring $R_2$ with finite weak
global dimension. Then, the direct product $R_1\times R_2$ is
GF-closed (by \cite[Proposition 3.5]{B}) with finite weak
Gorenstein global dimension (by Proposition \ref{prop-product}),
but it is neither coherent nor of finite weak global dimension.}
\end{exmp}

\end{section}



\bigskip\bigskip

\end{document}